\documentclass[12pt]{article}
\usepackage[T1]{fontenc}
\usepackage {lettrine}
\usepackage[portuges,english]{babel}
\usepackage{amsmath,latexsym,makeidx,amssymb,graphics,txfonts,color}
\usepackage[pdftex]{graphicx}
\usepackage[colorlinks,linkcolor=blue,urlcolor=blue,citecolor=blue,
plainpages=false,pdfpagelabels,breaklinks]{hyperref}

\definecolor{dk1}{rgb}{0.1,0.3,0.7}
\definecolor{dk2}{rgb}{0.8,0.1,0.3}
\definecolor{dk3}{rgb}{0.2,0.8,0.3}
\definecolor{glfc}{rgb}{0.3,0.1,0.9}
\definecolor{dk4}{rgb}{0.7,0.4,0.3}
\definecolor{gre}{rgb}{0.1,0.1,0.1}
\definecolor{blank}{rgb}{1.0,1.0,1.0}
\definecolor{ye1}{rgb}{0.6,0.8,0.1}

\title{\sc  A note on \\ Quasi-Ehresmann-Dedecker Universes}
\author{{\sc Décio Krause} \\ {\small Department of Philosophy} \\ {\small Federal University of Santa Catarina}}

\begin{document}
\maketitle

\newcommand{\igual}{:=}
\newcommand{\ita}{\textit}
\newtheorem{thm}{Theorem}[section]
\newtheorem{lem}{Lemma}[section]
\newtheorem{cor}{Corollary}[section]
\newtheorem{dfn}{Definition}[section]
\newtheorem{exe}{Example}[section]
\newtheorem{axm}{Axiom}[section]
\newcommand{\qst}{$\mathfrak{Q}$}
\newcommand{\cqd}{{\color{magenta}{\rule{.60ex}{1.7ex}}}}
\newcommand{\lra}{\leftrightarrow}
\newcommand{\Proof}{\noindent\textit{Proof:} \,}
\newcommand{\lan}{$\mathcal{L}_{\mathfrak{Q}}$}

\renewenvironment{enumerate}{\begin{list}{}{\rm \labelwidth 0mm
\leftmargin 5mm}} {\end{list}}

\newcommand{\qU}{\mathcal{U}}

Intuitively speaking, the purpose of employing universes in the foundations of set theory is that we can perform set-theoretical operations widely and still remain inside the universe. As is known, we cannot to it in theories such as ZFC, as Russell's paradox clearly shows. It may be thought that we could use classes as well, say in a theory like NBG, but yet in this case categories such as {\bf Grp} (groups), {\bf Ring}, {\bf Set} and other would be proper classes, the problem is that classes cannot be members of other classes, so no much gain would be obtained. The use of universes has been considered one of the best options (another would be the assumption of the existence of strongly inaccessible cardinals --- but this turns to be equivalent to the use of universes). 

We employ the resources of theory of quasi-sets $\mathfrak{Q}$ \cite{frekra10} form just sketching the category {\bf QSet} in a quite similar way that the category {\bf Set} is obtained from, say, the ZFC set-theory  by adding universes to it \cite{gro64}, \cite{bou64}. Since further developments are available for this study, in the sense that we may be interested in considering an universe as a starting point for other large universes (see \cite[pp. 259 and 262]{kro07} for the justification of Grotendieck's use of universes instead of NBG), we have opted for strengthening \qst\ with special universes, the \ita{quasi-Ehresmann-Dedecker universes}. 

Ehresmann-Dedecker universes were introduced in \cite{coscar67} (but see also \cite{cos67}, \cite{cos72}, \cite{car69}, \cite{bricos71}, \cite{car72}) and generalize the concept of Sonner-Grothendieck universes \cite{gro64}, \cite{bou64}, being more adapted for dealing with \ita{Urelemente}, which is the case of \qst.
We could exclude the axiom of regularity, but here we shall keep with it. We base much of our definitions and theorems in \cite{bricos71} and \cite{car72}. 

\begin{dfn}[qED universe] A non-empty quasi-set $\mathcal{U}$ is a \textit{quasi-Ehresmann-Dedecker universe} (qED for short) if the following conditions are obeyed:
\begin{enumerate}
\item (1) $x \in \qU \to \mathcal{P}(x) \in \qU$
\item (2) $x \in \qU \to   [x]_{\mathcal{U}} \in \qU$  
\item (3) $x, y \in \qU \to x \times y \in \qU$
\item (4) If $(x_i)_{i \in I}$ is a family of quasi-sets such that $x_i \in \qU$ for every $i \in I$ and being $I$ a "classical"\ quasi-set (that is, a quasi-set obeying the predicate $Z$), then $\bigcup_{i \in I} x_i \in \qU$.
\end{enumerate}
\end{dfn}

In (2),   $[x]_{\mathcal{U}}$     
stands for the quasi-set of all indistinguishable from $x$ that belong to the universe. It is called the \textit{singleton} of $x$. It should be remarked that this `singleton' may contain more than one element, that is, its cardinal may be greater that one. Intuitively speaking, indistinguishable objects cannot be distinguished by any means, but in \qst\ we may form collections (quasi-sets) of them with cardinals different from 1. In other words, indistinguishability doesn't collapse in identity, as in `Leibnizian' theories such as ZFC and others grounded on classical logic.

\begin{thm} If $\qU$ is a qED universe and $x, y \in \qU$, then $x \cup y$, $[x,y]_{\qU}$, and $\langle x,y \rangle_{\qU}$ belong to $\qU$. \\
\Proof {\rm The proof is similar, but with adaptations for the use of quasi-sets, to those presented in \cite[Teo.58, p.42]{car72}.} \cqd
\end{thm}

\textbf{Axiom AqED ---} Every quasi-set is an element of a qED universe.

\begin{dfn} $\mathfrak{Q}^\star := \mathfrak{Q} + AqED$ \end{dfn}

In words, $\mathfrak{Q}^\star$ is the theory got from  $\mathfrak{Q}$ by adding the axiom  {\bf AqED} to it. It is similar to the theory got from adding Grothendieck's universe axiom "For each set $x$, there exists a universe $\mathcal{U}$ such that $x \in \mathcal{U}$". 

\begin{thm} $\mathfrak{Q}^\star \vdash {\rm Cons}(\mathfrak{Q})$ \\
\Proof {\rm It is a homework to check that all postulates of $\mathfrak{Q}$ are satisfied in any qED $\qU$.} \cqd
\end{thm}

\begin{dfn} Let $\qU$ be a qED. We call $\qU$-{\bf qset} any element of $\qU$. A $\qU$-{\bf qclass} is a subset of $\qU$. A $\qU$-{\bf proper qclass} is a $\qU$-qclass that is not an $\qU$-qset. \end{dfn}

\begin{dfn} A $\qU$-{\bf small category} is a category $\mathbb{C}$ such that {\sf ob}$(\mathbb{C})$ and {\sf mor}$(\mathbb{C})$ are $\qU$-qsets, and it is $\qU$-{\bf large} otherwise. \end{dfn}

\begin{dfn} We call {\bf QSet} is the category of all quasi-sets of \qst. \end{dfn}

The objects of {\bf QSet} are the quasi-sets and the morphisms are the \ita{quasi-functions} between quasi-sets. We recall once more that the elements of a quasi-set may be indistinguishable from one each other and even so the cardinal of the collection (termed its \ita{quasi-cardinal}) is not identical to 1, as would be the case in `standard' set theories such as ZFC, NBG, KM or other, where the standard theory of identity implies that there cannot exist indistinguishable but not identical objects (these theories are, in a sense, \ita{Leibnizian}). Thus, taking two quasi-sets containing sub-collections of indistinguishable objects, a quasi-function takes indistinguishable objects in the first and associate to them indistinguishable elements of the other. In symbols, being $A$ and $B$ the quasi-sets and $q$ the quasi-function, we have that
$$(\forall x, x' \in A)(\exists y, y' \in B)(\langle x,y \rangle \in q \wedge \langle x', y' \rangle \in q \wedge x' \equiv x \to y \equiv y'),$$

\noindent where $\equiv$ is the relation of indistinguishability, which applies also to quasi-sets and in particular to quasi-functions. 

We can prove that the composition $\circ$ of quasi-functions (morphisms) is associative and that for each object $A$ there exists an `identity' morphism $A \stackrel{1_A}{\longrightarrow}A$, that is, a quasi-function whose domain and co-domain are both $A$ itself, such that for any other morphisms $B  \stackrel{f}{\longrightarrow}A$ and $A \stackrel{h}{\longrightarrow}C$, we have that $1_A \circ f \equiv f$ and $h \circ 1_A \equiv h$.  This characterizes the category {\bf Qset}. 

Since the collection of all quasi-sets and the collection of all quasi-functions are not quasi-sets, but things similar to proper classes, we can state the following result:

\begin{thm} The category {\bf QSet} is a large category. \end{thm}

The theorem shows that {\bf QSet} plays, relatively to quasi-set theory \qst\ the role played by {\bf Set} relatively to, say, ZFC. In other words, in \qst$^\star$  we can develop the category {\bf QSet} as suggested in \cite{dar14}.

\end{document}